\documentclass[11pt]{article}

\usepackage[utf8]{inputenc}
\usepackage{amsmath}
\usepackage{amsfonts,amssymb,latexsym,multicol,color}
\usepackage[francais]{babel}
\usepackage[T1]{fontenc}

\newcounter{fig}
\newtheorem{theo}{Th\'eor\`eme}

\newtheorem{prop}{Proposition}

\newcommand{\cad}{\text{c'est-\`a-dire }}
\newcommand{\expli}[1]{\quad\text{\footnotesize (#1)}}

\makeatletter
\ifcase\@ptsize

\or

\or

\fi
\makeatother

\newcommand{\Implique}{\Longrightarrow}

\newcommand{\fhi}{\varphi}
\newcommand{\ioe}{\leqslant}
\newcommand{\soe}{\geqslant}

\newcommand{\vers}{\rightarrow}

\newcommand{\demi}{{\frac{1}{2}}}

\newcommand{\ve}[1]{\boldsymbol{#1}}

\newcommand{\fin}{\hfill$\Box$}
\newcommand{\dem}{\noindent {\bf D\'emonstration\ }}
\newcommand{\fine}{\tag*{\mbox{$\Box$}}}

\providecommand{\bysame}{\leavevmode ---\ }
\providecommand{\og}{``}
\providecommand{\fg}{''}
\providecommand{\smfandname}{et}
\providecommand{\smfedsname}{\'eds.}
\providecommand{\smfedname}{\'ed.}
\providecommand{\smfmastersthesisname}{M\'emoire}
\providecommand{\smfphdthesisname}{Th\`ese}

\title{Sur la variation totale de la suite des parties fractionnaires des quotients d'un nombre r\'eel positif par les nombres entiers naturels cons\'ecutifs}

\author{Michel Balazard}
\date{}

\begin{document}
\maketitle

\begin{center}
  {\sc Abstract}
\end{center}
\begin{quote}
{\footnotesize We give an asymptotic formula for the total variation of the sequence of fractional parts of the quotients of a positive real number by the consecutive natural numbers :
$$
\sum_{n\soe 1}\lvert \{x/(n+1)\}-\{x/n\}\rvert=\frac{2}{\pi}\zeta(3/2)x^{1/2}+O(x^{2/5}).
$$}
\end{quote}

\begin{center}
  {\sc Keywords}
\end{center}
\begin{quote}
{\footnotesize Arithmetic functions, fractional part, total variation\\MSC classification : 11N37}
\end{quote}



\section{Introduction}

La quantit\'e d\'ecrite par le titre de cet article est
$$
W(x)=\sum_{n\soe 1}\lvert \{x/(n+1)\}-\{x/n\}\rvert \quad (x>0),
$$
o\`u
$\{t\}=t-\lfloor t\rfloor$ d\'esigne la partie fractionnaire du nombre r\'eel $t$, et $\lfloor t\rfloor$ sa partie enti\`ere. La lettre $W$ est choisie en r\'ef\'erence au math\'ematicien Aurel Wintner (1903-1958). Dans un article de 1946, \emph{Square root estimates of arithmetical sum functions} (cf. \cite{MR0016389}), il consid\'era la fonction $W(x)$ dans le contexte suivant.

\smallskip

Soit $f$ une fonction arithm\'etique \`a valeurs r\'eelles ou complexes dont la fonction sommatoire
\begin{equation}\label{t1}
F(x)=\sum_{n\ioe x}f(n) \quad (x>0)
\end{equation}
est born\'ee. En particulier, la s\'erie
$$
C=\sum_{n\soe 1}\frac{f(n)}{n}=\sum_{k\soe1}\frac{F(k)}{k(k+1)}=\int_0^{\infty}F(t)\frac{dt}{t^2}
$$
est alors convergente.

En notant $*$ le produit de convolution de Dirichlet des fonctions arithm\'etiques, posons $g=f*\ve{1}$ :
$$
g(n)=\sum_{d \mid n}f(d)\quad (n=1,2,\dots).
$$
On a
\begin{align*}
G(x)&=\sum_{n\ioe x}g(n)\\
&=\sum_{n\soe 1}f(n)\lfloor x/n\rfloor\\
&=Cx-\sum_{n\soe 1}f(n)\{x/n\}\\
&=Cx +\sum_{n\soe 1}F(n)\big (\{x/(n+1)\}-\{x/n\}\big).
\end{align*}

On en d\'eduit l'estimation
$$
\lvert G(x)-Cx\rvert\ioe \|F\|_{\infty}W(x) \quad (x>0),
$$
o\`u $\|F\|_{\infty}$ d\'esigne la borne sup\'erieure des modules des sommes \eqref{t1}. Dans \cite{MR0016389}, Wintner d\'emontra que 
\begin{equation}\label{t2}
W(x) \asymp x^{1/2} \quad (x\soe 1).
\end{equation}

On a donc $G(x)=Cx+O(x^{1/2})$. En outre, ce r\'esultat est optimal au sens o\`u, quelle que soit la fonction positive $\omega(x)$, d\'efinie pour $x>0$ et telle que
$$
\omega(x)=o(x^{1/2}) \quad (x\vers \infty),
$$
il existe une fonction arithm\'etique $f$, dont la fonction sommatoire est born\'ee, et pour laquelle la relation
\begin{equation}\label{t48}
G(x)=Cx+O(\omega(x))
\end{equation}
est fausse. Ce dernier fait d\'ecoule de la minoration $W(x)\gg x^{1/2}$ et du principe de condensation des singularit\'es (th\'eor\`eme de Banach et Steinhaus, cf. \cite{zbMATH02581680}\footnote{C'est le \emph{uniform boundedness principle} de la litt\'erature en langue anglaise ; Wintner invoqua le {\og Lebesgue-Toeplitz norm principle\fg}.}).

De plus Wintner remarqua que ces r\'esultats se d\'eduisaient uniquement de \eqref{t2}, et non de l'existence de la limite $\lim_{x\vers\infty} x^{-1/2}W(x)$, probl\`eme qu'il laissa ouvert et qui est l'origine du pr\'esent travail. 

\smallskip

Cela \'etant, on peut se passer complètement de l'introduction de la fonction $W$ et de considérations d'analyse fonctionnelle, aussi bien pour l'estimation $G(x)=Cx+O(x^{1/2})$ que pour l'optimalité de son terme d'erreur. D'une part, le principe de l'hyperbole de Dirichlet donne
$$
G(x)=\sum_{n\ioe \sqrt{x}}f(n)\lfloor x/n\rfloor +\sum_{m\ioe \sqrt{x}}F(x/m)-F(\sqrt{x})\lfloor \sqrt{x}\rfloor ,
$$
d'où découle simplement l'estimation de $G(x)$. D'autre part, Bayart a construit, dans le contexte de l'étude de l'abscisse de convergence d'un produit de séries de Dirichlet (cf. \cite{MR2039418}, \S2), un exemple d'une fonction arithm\'etique $f$ dont la fonction sommatoire est born\'ee, et pour laquelle l'assertion
$$
G(x)=Cx+o(x^{1/2}) \quad (x\vers \infty)
$$
est fausse ; la même fonction $f$ permet donc d'infirmer d'un coup toutes les assertions \eqref{t48}. Cette construction est rappelée au \S\ref{t37}.

\smallskip

L'argumentation de Wintner a cependant, entre autres mérites, celui d'attirer l'attention sur la question de l'estimation asymptotique de $W(x)$. 
\begin{theo}
Pour $x>0$, on a
$$
W(x)=\frac{2}{\pi}\zeta(3/2)\sqrt{x}+O(x^{2/5}),
$$
o\`u $\zeta$ d\'esigne la fonction $\zeta$ de Riemann.
\end{theo}

\smallskip

La nature de la somme $W(x)$ incite à l'analyse de la r\'epartition conjointe des deux fonctions
$$
n \mapsto \lfloor x/n\rfloor \quad \text{et} \quad n \mapsto \lfloor x/(n+1)\rfloor.
$$
Nous d\'eveloppons au \S\ref{t39} cette approche. Elle nous conduit \`a regrouper les termes de $W(x)$ suivant les valeurs $d$ de la diff\'erence $\lfloor x/n\rfloor-\lfloor x/(n+1)\rfloor$ ; le r\'esultat des calculs qui s'ensuivent est l'objet des propositions \ref{t6} et \ref{t35} (cf. \S\ref{t22} et \S\ref{t38} ci-dessous). Le th\'eor\`eme s'en d\'eduit par sommation au \S\S\ref{t40}-\ref{t41}.

La méthode employée est susceptible d'autres applications. Ainsi, des calculs similaires
fournissent l'énoncé suivant :
$$
\sum_{n\soe 1}\big (\{x/(n+1)\}-\{x/n\}\big)^2 = \frac{\zeta(3/2)}{\pi}\sqrt{x} +O(x^{3/7}) \quad (x>0).
$$

\section{D\'emonstration du th\'eor\`eme}\label{t39}

Dans toute la suite, la lettre $n$ désignera un nombre entier $\soe 1$ ; les lettres $h$, $k$, $d$ désigneront des nombres entiers $\soe 0$ ; la lettre $x$ désignera un nombre réel $>0$.

\subsection{R\'earrangement de la somme $W(x)$}

Nous r\'earrangeons la somme $W(x)$ suivant les valeurs de 
\begin{align*}
k&=\lfloor x/n\rfloor\\
h&=\lfloor x/(n+1)\rfloor
\end{align*}

On a donc $0\ioe h\ioe k \ioe x$ et
\begin{align*}
k\ioe x/n <k+1\\
h\ioe x/(n+1)<h+1,
\end{align*}
autrement dit
\begin{equation}\label{t3}
\max\big(\frac{x}{k+1},\frac{x}{h+1}-1\big)<n\ioe \min\big(\frac{x}{k},\frac{x}{h}-1\big).
\end{equation}
(avec la convention $x/0=\infty$). Nous d\'esignerons par $I(h,k;x)$ l'intervalle de valeurs de $n$ d\'efini par l'encadrement \eqref{t3}, \cad
$$
I(h,k;x)=\{1,2,\dots\}\,\cap\, ]x/(k+1),x/k]\,\cap\,]x/(h+1)-1,x/h-1].
$$
Il faut garder \`a l'esprit que, $x$ \'etant fix\'e, la collection des $I(h,k;x)$ non vides constitue une partition de l'ensemble des nombres entiers $\soe 1$.

Notons que, pour $n\in I(h,k;x)$, on a
$$
\{x/(n+1)\}-\{x/n\}= k-h-x/n(n+1).
$$
On a donc
$$
W(x)=\sum_{0\ioe h\ioe k\ioe x}W(h,k;x),
$$
o\`u
$$
W(h,k;x)=\sum_{n\in I(h,k;x)}\lvert k-h-x/n(n+1)\rvert.
$$

\smallskip

Maintenant, si $0\ioe d\ioe x$, nous posons
\begin{align}
E_d(x)&=\{(h,k) : 0\ioe h\ioe k\ioe x, \, k-h=d\}\label{t13}\\
&=\{(k-d,k) : d \ioe k \ioe x\}\notag
\end{align}
et
\begin{align*}
W_d(x)&=\sum_{(h,k)\in E_d(x)}W(h,k;x)\\
&=\sum_{d\ioe k\ioe x}\;\sum_{n\in I(k-d,k;x)}\lvert d-x/n(n+1)\rvert\, ,
\end{align*}
de sorte que
$$
W(x)=\sum_{0\ioe d\ioe x}W_d(x).
$$

Nous allons d'abord majorer la contribution \`a $W(x)$ des grandes valeurs de $d$, puis nous estimerons la quantit\'e $W_d(x)$, en commen\c{c}ant par le cas diagonal $d=0$.

\subsection{Contribution des grandes valeurs de $d$}

Soit $D$ un nombre r\'eel sup\'erieur \`a $1$. Si $d=\lfloor x/n\rfloor-\lfloor x/(n+1)\rfloor >D$, alors
\begin{align*}
\frac{x}{n^2} &>\frac{x}{n(n+1)}\\
&=d +\{ x/n\}-\{ x/(n+1)\}\\
&> D-1,
\end{align*}
donc $n<\sqrt{x/(D-1)}$. On en d\'eduit que
\begin{align}
\sum_{d>D}W_d(x)&\ioe\sum_{n<\sqrt{x/(D-1)}}1\notag\\
&<\sqrt{x/(D-1)}.\label{t12}
\end{align}

\subsection{Estimation de $W_0(x)$}\label{t22}

Nous avons 
$$
E_0(x)=\{(k,k) : 0\ioe k\ioe x\}.
$$
L'intervalle $I(k,k;x)$ est d\'efini par l'encadrement
\begin{equation}\label{t7}
\frac{x}{k+1}<n\ioe \frac{x}{k}-1.
\end{equation}
Par cons\'equent
\begin{align*}
W_0(x) &=\sum_{0\ioe k\ioe  x}W(k,k;x)\\
&=\sum_{0\ioe k\ioe  x}\quad\sum_{x/(k+1)<n\ioe x/k-1}x/n(n+1).
\end{align*}

Si $k(k+1) > x$, la somme int\'erieure est vide. D\'esignons donc par $K= K(x)$ le plus grand nombre entier $k$ tel que $k(k+1)\ioe x$, \cad
$$
K(x)=\lfloor\sqrt{x+1/4}-1/2\rfloor.
$$

Pour $x>0$, on a
\begin{align}
W_0 (x)&=\sum_{0\ioe k\ioe K}\quad\sum_{x/(k+1)<n\ioe x/k-1}x/n(n+1)\notag\\
&=x\sum_{0\ioe k\ioe K}\quad\sum_{x/(k+1)<n\ioe x/k}\frac 1{n(n+1)}-x\sum_{1\ioe k\ioe K}\frac{1}{\lfloor x/k \rfloor(\lfloor x/k \rfloor +1)}\notag\\
&=x\sum_{n>x/(K+1)}\frac{1}{n(n+1)}-x\sum_{1\ioe k\ioe K}\fhi(x/k),\label{t49}
\end{align}
o\`u l'on a pos\'e
\begin{align}
\fhi(t)&=\frac{1}{\lfloor t \rfloor(\lfloor t \rfloor+1)}\label{t19}\\
&=t^{-2}+O(t^{-3}) \quad (t\soe 1).\notag
\end{align}

\smallskip

D'une part, 
\begin{align}
x\sum_{n>x/(K+1)}\frac{1}{n(n+1)}&=\frac{x}{\lfloor x/(K+1) \rfloor+1}\notag\\
&=\frac{K+1}{1+(1-\{x/(K+1)\})(K+1)/x}\notag\\
&=K+O(1)\notag\\
&=\sqrt{x}+O(1).\label{t5}
\end{align}

D'autre part,
\begin{align}
x\sum_{1\ioe k\ioe K}\fhi(x/k)&=x\sum_{1\ioe k\ioe K}\big(k^2/x^2+O(k^3/x^3)\big)\notag\\
&=\frac{K^3+O(K^2)}{3x}+O(K^4/x^2)\notag\\
&=\frac 13\sqrt{x}+O(1).\label{t4}
\end{align}

\smallskip

La conjonction de \eqref{t49}, \eqref{t5} et \eqref{t4} donne le résultat suivant.
\begin{prop}\label{t6}
Pour $x>0$, on a
\begin{equation}
W_0(x)=\frac 23\sqrt{x}+O(1).
\end{equation}
\end{prop}

\subsection{Estimation de $W_d(x)$, $d>0$}

\subsubsection{D\'ecomposition de l'ensemble $E_d(x)$}

Si $d$ est positif, nous allons d\'ecomposer l'ensemble $E_d(x)$ d\'efini par \eqref{t13}
en une partition de trois sous-ensembles sur lesquels l'encadrement \eqref{t3} s'exprimera sans recours aux fonctions $\max$ et $\min$.

\smallskip

Si $k>h\soe 0$ et $x>0$, l'in\'egalit\'e
$$
\frac xk \ioe \frac xh -1
$$
\'equivaut \`a
$$
\frac{hk}{k-h} \ioe x.
$$
En particulier, on  a les implications
$$
\frac{x}{k+1}\ioe\frac{x}{h+1}-1 \Implique \frac xk \ioe \frac xh -1
$$
et
$$
\frac xk > \frac xh -1 \Implique \frac{x}{k+1} > \frac{x}{h+1}-1.
$$

Cela nous incite \`a consid\'erer les trois parties suivantes de $E_d(x)$  (la d\'efinition de chaque $E_{d,i}(x)$ est suivie par la forme que prend l'encadrement \eqref{t3} lorsque $(k-d,k)\in E_{d,i}(x)$) :
\begin{align}
E_{d,1}(x)&=\{(k-d,k) : d\ioe k\ioe x, \, (k-d+1)(k+1) \ioe dx\}\notag\\
&\quad \quad\quad \frac{x}{k-d+1}-1<n\ioe \frac{x}{k}\label{t8}\\ \notag\\
E_{d,2}(x)&=\{(k-d,k) : d\ioe k\ioe x, \, (k-d)k > dx\}\notag\\
&\quad\quad\quad\frac{x}{k+1}<n\ioe \frac{x}{k-d}-1\label{t9}\\ \notag\\
E_{d,3}(x)&=\{(k-d,k) : d\ioe k\ioe x, \, (k-d)k\ioe dx <(k-d+1)(k+1)\}\notag\\
&\quad\quad\quad\frac{x}{k+1}<n\ioe \frac{x}{k}\label{t10}
\end{align}

\smallskip

Celles des trois parties $E_{d,i}(x)$ ($1\ioe i\ioe 3$) qui sont non vides forment une partition de $E_d(x)$. Par cons\'equent, on a
$$
W_d(x)=W_{d,1}(x)+W_{d,2}(x)+W_{d,3}(x),
$$
o\`u
$$
W_{d,i}(x)=\sum_{(h,k)\in E_{d,i}(x)} W(h,k;x) \quad (1\ioe i \ioe 3).
$$

Avant d'\'evaluer successivement les trois quantit\'es $W_{d,i}(x)$, nous allons aux paragraphes suivants d\'efinir et \'etudier deux fonctions auxiliaires, $K_d(x)$ et $N_d(x)$.

\subsubsection{La fonction $K_d(x)$}

Pour $x > 0$ et $d\soe 0$, nous d\'efinissons $K_d(x)$ comme le plus grand nombre entier $k$ tel que
$$
(k-d)k\ioe dx\, ,
$$
in\'egalit\'e qui \'equivaut \`a
$$
\frac xk \ioe \frac{x}{k-d} -1
$$
si $k>d$.

On a donc $K_0(x)=0$ et en général
\begin{equation*}
K_d(x)=\lfloor(d+\sqrt{d^2+4dx})/2\rfloor.
\end{equation*}

En utilisant le fait que $t \mapsto (t-d)t$ est strictement croissante sur $[d/2,\infty[$, on démontre les relations 
\begin{align}
d \ioe K_d(x)&\ioe x+d  \notag\\ 
d+1 \ioe K_d(x)&  \quad\qquad\qquad(d>0, \, x\soe 2)\notag\\
2d\ioe K_d(x)&\ioe x  \quad\qquad (x\soe 2d)\label{t50}\\
\sqrt{dx}+d/2-1< K_d(x)&\ioe \sqrt{dx}+d \label{t23}\\
K_d(x)&<K_{d+1}(x).\notag
\end{align}

\smallskip

Nous utiliserons de plus des estimations des sommes $\sum 1$ et $\sum k^2$ portant sur les nombres entiers $k$ de l'intervalle $]K_{d}(x),K_{d+1}(x)]$.
\begin{prop}\label{t29}
Pour $0\ioe d\ioe x$ et $x\soe 1$, on a
\begin{equation*}
K_{d+1}(x)-K_{d}(x)=\big(\sqrt{d+1}-\sqrt{d}\, \big)\sqrt{x}+O(1).
\end{equation*}
\end{prop}
\dem

Comme fonction de $d$, la quantit\'e 
$$
\sqrt{d^2+4dx}-\sqrt{4dx}=\frac{1}{\sqrt{1/d^2+4x/d^3}+\sqrt{4x/d^3}}
$$
est croissante. D'autre part, sa d\'eriv\'ee par rapport \`a $d$ est
\begin{align}
 \frac{d+2x}{\sqrt{d^2+4dx}}-\sqrt{\frac xd}&=\frac{d\sqrt{d}+2x\sqrt{d}-\sqrt{x(d^2+4dx)}}{\sqrt{d(d^2+4dx)}}\notag\\
&\ioe \demi\sqrt{\frac dx}.\label{t28}
\end{align}

Par cons\'equent
\begin{align*}
K_{d+1}(x)-K_d(x)&=\big(d+1+\sqrt{(d+1)^2+4(d+1)x}\,\big)/2-\big(d+\sqrt{d^2+4dx}\,\big)/2 +O(1)\\
&=\demi\big(\sqrt{4(d+1)x}-\sqrt{4dx}\, \big)+O\big(\sqrt{(d+1)/x}\, \big) +O(1)\expli{d'apr\`es \eqref{t28}}\\
&=\big(\sqrt{d+1}-\sqrt{d}\, \big)\sqrt{x}+O(1),
\end{align*}
si $d\ioe x$ et $x\soe 1$.\fin

\smallskip

Notons que la proposition \ref{t29} entraîne l'estimation
$$
K_{d+1}(x)-K_d(x)\ll \sqrt{x/(d+1)} \quad (0\ioe d\ioe x, \, x\soe 1).
$$

\begin{prop}\label{t51}
Pour $0\ioe d\ioe x$ et $x\soe 1$, on a
\begin{equation*}
\sum_{K_{d}(x)<k\ioe K_{d+1}(x)}k^2=\frac{(d+1)\sqrt{d+1}-d\sqrt{d}}{3}\,x^{3/2}+O\big((d+1)x\big).
\end{equation*}
\end{prop}
\dem

Si $M$ et $N$ sont des entiers naturels tels que $M\ioe N$, on a
$$
\sum_{M<k\ioe N}k^2=\frac{N-M}{6}\big(2(N^2+NM+M^2)+3(N+M)+1\big).
$$
Avec $M=K_{d}(x)$ et $N=K_{d+1}(x)$, et compte tenu de \eqref{t23} et de la proposition \ref{t29}, cela donne, 
\begin{align*}
\sum_{K_{d}<k\ioe K_{d+1}}k^2&=\frac{(\sqrt{d+1}-\sqrt{d})\sqrt{x}+O(1)}{6}\big(2(2d +1+\sqrt{d(d+1)})x+O\big((d+1)^{3/2}\sqrt{x})\big)\\
&=\frac{(d+1)\sqrt{d+1}-d\sqrt{d}}{3}\,x^{3/2}+O\big((d+1)x\big).\fine
\end{align*}

\smallskip

Au moyen de la fonction $K_d(x)$, on peut r\'ecrire les conditions, quadratiques relativement \`a $k$, intervenant dans les d\'efinitions des ensembles $E_{d,i}(x)$, sous les formes suivantes, respectivement :
\begin{align}
k&\ioe K_d(x)-1 & (i=1)\label{t16}\\
k&>K_d(x) & (i=2)\label{t17}\\
k&=K_d(x) & (i=3)\label{t18}
\end{align}

\subsubsection{La fonction $N_d(x)$}

Pour exprimer la quantit\'e $\lvert d-x/n(n+1)\rvert$ sans valeur absolue, nous sommes conduits \`a d\'efinir, pour $d>0$, $N_d(x)$ comme le plus grand nombre entier tel que
$$
n(n+1)\ioe x/d.
$$
On a donc
$$
N_d(x)=\lfloor(-1+\sqrt{1+4x/d}\,)/2\rfloor.
$$

Notons les relations suivantes.
\begin{align}
 1&\ioe N_d(x)  &(0 <d\ioe x/2)\notag\\ 
N_d(x)&\ioe \sqrt{x/d}  &(d>0, x\soe 0)\label{t32}\\
\sqrt{x/d}-2 &< N_d(x)   &(0<d\ioe x)\label{t33}
\end{align}

\medskip

\'Etablissons maintenant une relation entre les fonctions $K_d$ et $N_d$ (pour ces deux fonctions, nous omettons dor\'enavant la mention de la variable $x$ afin d'all\'eger les notations).
\begin{prop}\label{t14}
Pour $d$ entier et $x$ r\'eel tels que $0<d\ioe x$, on a
$$
\lfloor x/(K_d+1)\rfloor \ioe N_d \ioe \lfloor x/K_d\rfloor.
$$
\end{prop}
\dem

Par d\'efinition de $N_d$, et comme $t\mapsto t(t+1)$ est strictement croissante pour $t\soe 0$, il suffit de v\'erifier que
$$
\frac{x}{K_d+1}\Big( \frac{x}{K_d+1} +1\Big) < \frac xd \ioe \frac{x}{K_d}\Big( \frac{x}{K_d} +1\Big).
$$
Or cet encadrement \'equivaut au suivant :
$$
(K_d-d)K_d \ioe dx < (K_d+1-d)(K_d+1)\, ,
$$
lequel d\'ecoule de la d\'efinition de $K_d$.\fin

\subsubsection{Calcul de $W_{d,1}(x)$}

En supposant $0<d\ioe x/2$, on a d'apr\`es \eqref{t8}, \eqref{t50} et \eqref{t16} :
\begin{equation*}
W_{d,1}(x)=\sum_{d\ioe k\ioe K_d -1}\;\sum_{\frac{x}{k-d+1}-1<n\ioe \frac{x}{k}}\lvert d-x/n(n+1)\rvert.
\end{equation*}

La somme int\'erieure est non vide seulement si 
$$
\frac{x}{k-d+1}-1< \frac{x}{k}\, ,
$$
autrement dit seulement si $k>K_{d-1}$ (rappelons que $K_0=0$). Les nombres entiers $n$ intervenant dans cette somme int\'erieure sont strictement sup\'erieurs \`a
$$
\frac{x}{K_d-d}-1 \soe \frac{x}{K_d} \soe N_d,
$$
d'apr\`es la d\'efinition de $K_d$ et la proposition \ref{t14}. 

Dans le calcul qui suit, ainsi qu'au paragraphe suivant, nous \'ecrirons comme au \S\ref{t22} :
$$
\frac{1}{\lfloor t\rfloor +1}=\frac{1}{\lfloor t\rfloor}-\fhi(t)\, ,
$$
et nous emploierons l'identit\'e {\og $r$-t\'elescopique\fg} :
$$
\sum_{a<k\ioe b} (u_k-u_{k-r})= \sum_{b-r<k\ioe b} u_k-\sum_{a-r<k\ioe a} u_k.
$$

\smallskip

On a donc
\begin{multline}\label{t20}
W_{d,1}(x)=\sum_{K_{d-1} < k\ioe K_d -1}\;\sum_{\frac{x}{k-d+1}-1<n\ioe \frac{x}{k}}\Big (d- \frac{x}{n(n+1)}\Big)\\
=\sum_{K_{d-1} < k\ioe K_d -1}\; \bigg(d\big(\lfloor x/k \rfloor-\lfloor x/(k-d+1)\rfloor+1\big)-x\Big(\frac{1}{\lfloor x/(k-d+1)\rfloor}-\frac{1}{\lfloor x/k\rfloor+1}\Big)\bigg)\\
=d\sum_{K_{d} -d< k\ioe K_d -1}\; \lfloor x/k \rfloor-d\sum_{K_{d-1} -d+1< k\ioe K_{d -1}}\; \lfloor x/k \rfloor+d(K_d-K_{d-1}-1)+\\
+x\sum_{K_{d} -d< k\ioe K_d -1}\; \frac{1}{\lfloor x/k \rfloor}-x\sum_{K_{d-1} -d+1< k\ioe K_{d -1}}\; \frac{1}{\lfloor x/k \rfloor}- x\sum_{K_{d-1} < k\ioe K_d -1}\;\fhi(x/k)
\end{multline}
o\`u $\fhi$ est d\'efinie par \eqref{t19}. Observons que, dans le cas o\`u $d=1$, les quatre premi\`eres sommes $\sum$ du résultat de \eqref{t20} sont vides.

\subsubsection{Calcul de $W_{d,2}(x)$}

En supposant toujours $0<d\ioe x/2$, on a d'apr\`es \eqref{t9} et \eqref{t17} :
\begin{equation*}
W_{d,2}(x) =\sum_{K_d <k \ioe x}\;\sum_{\frac{x}{k+1}<n\ioe \frac{x}{k-d}-1}\lvert d-x/n(n+1)\rvert
\end{equation*}

La somme int\'erieure est non vide seulement si 
$$
\frac{x}{k-d}-1\soe \frac{x}{k+1}\, ,
$$
autrement dit seulement si $k\ioe K_{d+1}-1$. Les nombres entiers $n$ intervenant dans cette somme int\'erieure sont inf\'erieurs ou \'egaux \`a
$$
\Big\lfloor\frac{x}{K_d+1-d}-1\Big\rfloor \ioe \Big\lfloor\frac{x}{K_d+1}\Big\rfloor \ioe N_d,
$$
d'apr\`es la d\'efinition de $K_d$ et la proposition \ref{t14}.

En supposant $d+1\ioe x/2$, on a donc
\begin{multline}\label{t21}
W_{d,2}(x) =\sum_{K_d <k \ioe K_{d+1}-1}\;\sum_{\frac{x}{k+1}<n\ioe \frac{x}{k-d}-1} \big(x/n(n+1)-d\big)\\
=\sum_{K_{d} < k\ioe K_{d+1} -1}\; \bigg(x\Big(\frac{1}{\lfloor x/(k+1)\rfloor +1}-\frac{1}{\lfloor x/(k-d)\rfloor}\Big)-d\big(\lfloor x/(k-d) \rfloor-\lfloor x/(k+1)\rfloor-1\big)\bigg)\\
=d\sum_{K_{d+1} -d-1 < k\ioe K_{d+1}}\; \lfloor x/k \rfloor  -d\sum_{K_{d} -d< k\ioe K_{d} +1}\; \lfloor x/k \rfloor +d(K_{d+1}-K_{d}-1)+\\
+x\sum_{K_{d+1} -d-1< k\ioe K_{d+1}}\; \frac{1}{\lfloor x/k \rfloor}-x\sum_{K_{d} -d< k\ioe K_{d}+1}\; \frac{1}{\lfloor x/k \rfloor} -x\sum_{K_{d} +1< k\ioe K_{d+1} }\;\fhi(x/k).\end{multline}

\subsubsection{Calcul de $W_{d,3}(x)$}

Pour $0<d\ioe x$, on a d'apr\`es \eqref{t10} et \eqref{t18} :
\begin{align*}
W_{d,3}(x)&=\sum_{\frac{x}{K_d+1}<n\ioe \frac{x}{K_d}}\lvert d-x/n(n+1)\rvert\\
&=W_{d,3}^-(x)+W_{d,3}^+(x),
\end{align*}
o\`u
\begin{align*}
W_{d,3}^-(x)&=\sum_{\frac{x}{K_d+1}<n\ioe N_d}\Big( \frac{x}{n(n+1)}-d\Big)\\
W_{d,3}^+(x)&=\sum_{N_d<n\ioe \frac{x}{K_d}}\Big(d- \frac{x}{n(n+1)}\Big)
\end{align*}
d'apr\`es la proposition \ref{t14}. La somme $W_{d,3}^-(x)$ est vide si $N_d=\lfloor x/(K_d+1)\rfloor$.

\smallskip

On a
\begin{align}
W_{d,3}^-(x)&=x\Big( \frac{1}{\lfloor x/(K_d+1)\rfloor +1}-\frac{1}{N_d+1}\Big) -d(N_d-\lfloor x/(K_d+1)\rfloor )\notag\\
W_{d,3}^+(x)&=d(\lfloor x/K_d\rfloor -N_d)-x\Big( \frac{1}{N_d+1}-\frac{1}{\lfloor x/K_d\rfloor +1}\Big)\notag\\
W_{d,3}(x)&=x\Big( \frac{1}{\lfloor x/(K_d+1)\rfloor +1}+\frac{1}{\lfloor x/K_d\rfloor +1}-\frac{2}{N_d+1}\Big) -d(2N_d-\lfloor x/(K_d+1)\rfloor -\lfloor x/K_d\rfloor).\label{t15}
\end{align}

\subsubsection{Calcul et estimation de $W_d(x)$}\label{t38}

En ajoutant \eqref{t20}, \eqref{t21} et \eqref{t15} et en r\'eduisant, on obtient pour $0<d\ioe x/2-1$ :
\begin{multline}\label{t26}
W_d(x) =W_{d,1}(x)+W_{d,2}(x)+W_{d,3}(x)\\
=\sum_{K_{d+1} -d-1 < k\ioe K_{d+1}}\big( d\lfloor x/k \rfloor +x/\lfloor x/k \rfloor\big)-\sum_{K_{d-1} -d+1< k\ioe K_{d -1}}\big(d\lfloor x/k \rfloor +x/\lfloor x/k \rfloor\big)+\\
-x\sum_{K_{d-1} < k\ioe K_{d+1} }\;\fhi(x/k) +d(K_{d+1}-K_{d-1}-2)-\frac{2x}{N_d+1} -2dN_d.
\end{multline}

\smallskip

Pour estimer les deux premières sommes de \eqref{t26}, nous utiliserons la proposition suivante.
\begin{prop}\label{t44}
Pour $0<d\ioe x/2$ et $K_d-d<k\ioe K_d$, on a
$$
d\lfloor x/k \rfloor +x/\lfloor x/k \rfloor=2\sqrt{dx}+O(d^{3/2}x^{-1/2}).
$$
\end{prop}
\dem

Posons $q=\lfloor x/k \rfloor$. On a
$$
dq+x/q-2\sqrt{dx}=p^2,
$$
où
\begin{align*}
p&= \sqrt{dq}-\sqrt{x/q}\\
&=\frac{dq^2-x}{q\big(\sqrt{dq}+\sqrt{x/q}\big)}.
\end{align*}

Or
\begin{align}
q&=\lfloor x/k \rfloor\notag\\
&=x/k+O(1)\notag\\
&=\frac{x}{\sqrt{dx}+O(d)}+O(1)\expli{d'après \eqref{t23}}\notag\\
&=\sqrt{x/d}+O(1).\label{t45}
\end{align}
Par conséquent,
\begin{align*}
dq^2-x&=O(\sqrt{dx})\\
p&=O(d^{3/4}x^{-1/4})
\end{align*}
et 
\begin{equation*}
dq+x/q-2\sqrt{dx}=O(d^{3/2}x^{-1/2}).\fine
\end{equation*}

\smallskip

Nous pouvons maintenant démontrer le résultat suivant.
\begin{prop}\label{t35}
Pour $d >0$ et $x>0$, on a
$$
W_d(x)= f(d)\sqrt{x}+O(d),
$$
o\`u
\begin{equation*}
f(d) =\frac{8d+2}3\sqrt{d+1}-\frac{8d-2}3\sqrt{d-1}-4\sqrt{d}
\end{equation*}
\end{prop}
\dem

Nous appliquons la proposition \ref{t44} en changeant $d$ en $d +1$ et obtenons, si $0<d\ioe x/2-1$, 
\begin{multline}\label{t25}
\sum_{K_{d+1} -d-1 < k\ioe K_{d+1}}\big( d\lfloor x/k \rfloor +x/\lfloor x/k \rfloor\big)=\\
-\sum_{K_{d+1} -d-1 < k\ioe K_{d+1}}\lfloor x/k \rfloor +\sum_{K_{d+1} -d-1 < k\ioe K_{d+1}}\big(2\sqrt{(d+1)x}+O(d^{3/2}x^{-1/2})\big)\\
=(2d+1)\sqrt{(d+1)x}+O(d)+O(d^{5/2}x^{-1/2}),
\end{multline}
où l'on a utilisé \eqref{t45} (avec $d$ remplacé par $d+1$) pour évaluer l'avant-dernière somme.

De même,
\begin{multline}\label{t27}
\sum_{K_{d-1} -d+1 < k\ioe K_{d-1}}\big( d\lfloor x/k \rfloor +x/\lfloor x/k \rfloor\big)=\\
\sum_{K_{d-1} -d+1 < k\ioe K_{d-1}}\lfloor x/k \rfloor +\sum_{K_{d-1} -d+1 < k\ioe K_{d-1}}\big(2\sqrt{(d-1)x}+O(d^{3/2}x^{-1/2})\big)\\
=(2d-1)\sqrt{(d-1)x}+O(d)+O(d^{5/2}x^{-1/2}),
\end{multline}
résultat valable même si $d=1$.

\smallskip

Maintenant, d'après la proposition \ref{t51}, on a
\begin{equation*}
\sum_{K_{d-1}<k\ioe K_{d+1}}k^2=\frac{(d+1)\sqrt{d+1}-(d-1)\sqrt{d-1}}{3}x^{3/2}+O(dx).
\end{equation*}
De plus,
\begin{align*}
\sum_{K_{d-1}<k\ioe K_{d+1}}k^3&\ll \sqrt{x/d}\cdot \sqrt{(dx)^3}\\
&=dx^2.
\end{align*}

Par cons\'equent,
\begin{align}
x\sum_{K_{d-1} < k\ioe K_{d+1}}\;\fhi(x/k)&=\frac 1x\sum_{K_{d-1} < k\ioe K_{d+1}}k^2+O(x^{-2}\sum_{K_{d-1} < k\ioe K_{d+1}}k^3)\notag\\
&=\frac{(d+1)\sqrt{d+1}-(d-1)\sqrt{d-1}}{3}\sqrt{x}+O(d).\label{t30}
\end{align}

\smallskip

Le terme suivant de \eqref{t26} peut être estimé grâce à la proposition \ref{t29}. On a
\begin{equation}\label{t31}
d(K_{d+1}-K_{d-1}-2)=(d\sqrt{d+1}-d\sqrt{d-1})\sqrt{x}+O(d).
\end{equation}

\smallskip

Enfin, en utilisant \eqref{t32} et \eqref{t33}, on a
\begin{align}
-\frac{2x}{N_d+1} -2dN_d&=-2\frac{x}{\sqrt{x/d}+O(1)}-2d\big(\sqrt{x/d}+O(1)\big)\notag\\
&=-4\sqrt{dx}+O(d)\label{t34}
\end{align}

\smallskip

En insérant maintenant \eqref{t25}, \eqref{t27}, \eqref{t30}, \eqref{t31} et \eqref{t34} dans \eqref{t26}, nous obtenons 
\begin{equation}\label{t46}
W_d(x)= f(d)\sqrt{x}+O(d) +O(d^{5/2}x^{-1/2}) \quad (0<d\ioe x/2-1).
\end{equation}

Si $d\ioe x^{1/3}$, le second terme d'erreur est absorbé par le premier. D'autre part, en utilisant l'approximation 
$$
\sqrt{d\pm 1}=\sqrt{d}\big(1\pm 1/2d-1/8d^2+O(d^{-3})\big)\, ,
$$
on voit que $f(d)=O(d^{-3/2})$. La majoration uniforme $W_d(x)\ll \sqrt{x/d}$ montre alors
que
$$
W_d(x)- f(d)\sqrt{x} \ll d \quad (d>x^{1/3}),
$$
ce qui permet d'omettre définitivement le second terme d'erreur de \eqref{t46} et la condition $d\ioe x/2-1$.\fin

\subsection{Sommation de la s\'erie des $f(d)$}\label{t40}

Puisque $f(d)=O(d^{-3/2})$, la s\'erie $\sum_{d\soe 1}f(d)$ converge ; nous allons calculer sa somme.

\smallskip

En \'ecrivant
$$
f(d)=\frac{8(d+1)-6}3\sqrt{d+1}-\frac{8(d-1)-6}3\sqrt{d-1}-4\big (\sqrt{d-1}+\sqrt{d}\,\big)\, ,
$$
et en supposant $D$ entier positif, on obtient
$$
\sum_{d=1}^D f(d)=\frac{8(D+1)-6}3\sqrt{D+1} + \frac{8D-6}3\sqrt{D}-\frac 23-8\sum_{d=1}^D\sqrt{d}+4\sqrt{D}.
$$

Or
$$
\frac{8(D+1)-6}3\sqrt{D+1}=\frac 83 D^{3/2}+2D^{1/2}+O(D^{-1/2})
$$
et la formule sommatoire d'Euler et Maclaurin donne
$$
\sum_{d=1}^D\sqrt{d}=\frac 23 D^{3/2}+\demi D^{1/2}+\zeta(-1/2)+O(D^{-1/2})
$$
(cf. \cite{hardyDS} (13$\cdot$10$\cdot$7), p. 333, et \cite{B} chapter 7, (1$\cdot$2), p. 150, (4$\cdot$5), p. 156). On en d\'eduit
\begin{align}
\sum_{d=1}^{\infty} f(d)&=-\frac 23-8\zeta(-1/2)\notag\\
&=-\frac 23 +\frac{2}{\pi}\zeta(3/2),\label{t36}
\end{align}
d'apr\`es l'\'equation fonctionnelle de la fonction $\zeta$ :
$$
\zeta(1-s)=2(2\pi)^{-s}\cos(\pi s/2)\Gamma(s)\zeta(s).
$$

\subsection{Conclusion}\label{t41}

Pour $0<x\ioe 1$, on a $W(x)=x$. Pour $x>1$ et $D\soe 2$, on a
\begin{align*}
W(x) &= W_0(x)+\sum_{1\ioe d\ioe D}W_d(x)+\sum_{d>D}W_d(x)\\
&=\frac 23\sqrt{x}+O(1)+\sum_{1\ioe d\ioe D}\big(f(d)\sqrt{x}+O(d)\big) +O(\sqrt{x/D})\expli{d'apr\`es \eqref{t12} et les propositions \ref{t6} et \ref{t35}}\\
&=\big (2/3+\sum_{d=1}^{\infty} f(d)\big)\sqrt{x}+O(D^2)+O\big(\sqrt{x/D}\,\big)\expli{car $f(d)=O(d^{-3/2})$}\\
&=\frac{2}{\pi}\zeta(3/2)\sqrt{x}+O(x^{2/5}),
\end{align*}
d'apr\`es \eqref{t36}, et en choisissant $D=2x^{1/5}$.

\section{Une remarque sur l'optimalit\'e du terme d'erreur $O(\sqrt{x})$}\label{t37}

Au \S2 de \cite{MR2039418}, Bayart d\'efinit une fonction arithm\'etique $\alpha$ de la fa\c{c}on suivante. On construit d'abord une suite de nombres entiers par blocs en posant
\begin{align*}
M_n &=2^{4^n} \quad (n\soe 1)\\
i_{k,n} &=\lfloor M_n/k\rfloor \quad (1\ioe k \ioe \sqrt{M_n}/2=2^{2^{2n-1}-1})
\end{align*} 

Le plus petit \'el\'ement $i_{k,n}$ du $n$\up{e} bloc correspond \`a $k=\sqrt{M_n}/2$, et vaut $2\sqrt{M_n}$. On montre que les nombres $i_{k,n}$ v\'erifient l'\'egalit\'e $\lfloor M_n/i_{k,n}\rfloor =k$ ; ils forment donc une suite strictement d\'ecroissante pour chaque valeur de $n$.

On pose ensuite
$$
\alpha(m)=\begin{cases}
(-1)^k &\text{ si existent $n$ et $k$ tels que $m=i_{k,n}$}\\
0 &\text{ sinon.}
\end{cases}
$$ 

Pour tout $x>0$ on a donc $A(x)=\sum_{n\ioe x}\alpha(n)\in \{0,1\}$. 

\smallskip

Consid\'erons maintenant la fonction arithm\'etique $\beta$ définie par $\beta(n)=(-1)^{n-1}$. De même que pour $\alpha$, la fonction sommatoire $B(x)$ de $\beta$ ne prend que les valeurs $0$ et $1$ ; pour $k$ entier, on a $B(k)=0$ si $k$ est pair, et $B(k)=1$ si $k$ est impair.

Posons ensuite $g=\alpha*\beta$. D'apr\`es le principe de l'hyperbole de Dirichlet, la fonction sommatoire $G$ de $g$ v\'erifie 
$$
G(x) = \sum_{n\ioe \sqrt{x}}\alpha(n)B(x/n)+\sum_{n\ioe \sqrt{x}}\beta(n)A(x/n)-B(\sqrt{x})A(\sqrt{x}),
$$
donc
\begin{equation}\label{t42}
\lvert G(x) \rvert \ioe 2\sqrt{x}+1 \quad (x>0)
\end{equation}

Or pour $x=M_n$,
\begin{align*}
G(x) &=\sum_{i\ioe x} \alpha(i)B(x/i)\\
 &=\sum_{i\ioe 2\sqrt{x}-1} +\sum_{2\sqrt{x}\ioe i \ioe x}\\
&=S_1+S_2, \text{ disons.}
\end{align*}

Dans $S_1$, tous les termes correspondants \`a $i>M_{n-1}=x^{1/4}$ sont nuls ; par cons\'equent $\vert S_1\rvert \ioe x^{1/4}$.

Dans $S_2$ les termes $\alpha(i)$ non nuls correspondent exactement aux \'el\'ements $i_{k,n}$, et on a alors $B(x/i_{k,n})=B(k)$. Par cons\'equent,
\begin{align*}
S_2 &=\sum_{\substack{
1\ioe k\ioe \sqrt{x}/2\\ 
k\text{ impair}
}
}(-1)^k\\
&=-\sqrt{x}/4.
\end{align*}

On a donc
\begin{equation}\label{t43}
G(M_n)=\sum_{i\ioe M_n} \alpha(i)B(M_n/i) \ioe -\sqrt{M_n}/4+M_n^{1/4} \quad (n\soe 1).
\end{equation}

\smallskip

Maintenant, observons que $\beta=\gamma*\ve{1}$, où $\gamma$ est la fonction arithmétique définie par $\gamma(1)=1$, $\gamma(2)=-2$ et $\gamma(n)=0$ pour $n>2$, et posons $f=\alpha*\gamma$, de sorte que $g=\alpha*\beta=f*\ve{1}$. 

Comme la fonction sommatoire de $f$, \`a savoir $F(x)=A(x)-2A(x/2)$, est born\'ee, nous sommes dans la situation d\'ecrite dans l'introduction de cet article. On a donc
\begin{equation}\label{t47}
G(x)=Cx +O(\sqrt{x}),
\end{equation}
avec $C=\sum_nf(n)/n$. La relation \eqref{t42} prouve que $C=0$, et \eqref{t43} prouve que le $O(\sqrt{x})$ de \eqref{t47} ne peut pas \^etre remplac\'e par un $o(\sqrt{x})$.



\begin{thebibliography}{1}

\bibitem{zbMATH02581680}
{\scshape S.~{Banach} {\normalfont \smfandname} H.~{Steinhaus}} -- {\og {Sur le
  principe de la condensation de singularit\'es.}\fg}, \emph{{Fundam. Math.}}
  \textbf{9} (1927), p.~50--61.

\bibitem{MR2039418}
{\scshape F.~Bayart} -- {\og The product of two {D}irichlet series\fg},
  \emph{Acta Arith.} \textbf{111} (2004), p.~141--152.

\bibitem{B}
{\scshape B.~C. Berndt} -- \emph{Ramanujan notebooks, part 1}, Springer,
  Berlin, 1985.

\bibitem{hardyDS}
{\scshape G.~H. Hardy} -- \emph{Divergent series}, Oxford University Press,
  1949.

\bibitem{MR0016389}
{\scshape A.~Wintner} -- {\og Square root estimates of arithmetical sum
  functions\fg}, \emph{Duke Math. J.} \textbf{13} (1946), p.~185--193.

\end{thebibliography}

\providecommand{\bysame}{\leavevmode ---\ }
\providecommand{\og}{``}
\providecommand{\fg}{''}
\providecommand{\smfandname}{et}
\providecommand{\smfedsname}{\'eds.}
\providecommand{\smfedname}{\'ed.}
\providecommand{\smfmastersthesisname}{M\'emoire}
\providecommand{\smfphdthesisname}{Th\`ese}

\medskip

\footnotesize

\noindent BALAZARD, Michel\\
Aix Marseille Université, CNRS, Centrale Marseille, I2M UMR 7373\\
13453, Marseille\\
FRANCE\\
Adresse \'electronique : \texttt{balazard@math.cnrs.fr}

\end{document}